\documentclass[twoside,leqno,10pt]{article}
\usepackage{graphics,ifthen}
\usepackage{latexsym}
\usepackage{amsmath}
\usepackage{amssymb}
\usepackage{mathrsfs}
\begin{document}
\input{latexP.sty}
\input{referencesP.sty}
\input epsf.sty
\def\ind{\stackrel{\mathrm{ind}}{\sim}}
\def\iid{\stackrel{\mathrm{iid}}{\sim}}
\def\Definition{\stepcounter{definitionN}\
    \Demo{Definition\hskip\smallindent\thedefinitionN}}
\def\EndDefinition{\EndDemo}
\def\Example#1{\Demo{Example [{\rm #1}]}}
\def\EndExample{\qed\EndDemo}
\def\Category#1{\centerline{\Heading #1}\rm}
\
%% Paper specific definitions
\def\e{\text{\hskip1.5pt e}}
\newcommand{\eps}{\epsilon}
\newcommand{\proof}{\noindent {\bf Proof:\ }}
\newcommand{\remarks}{\noindent {\bf Remarks:\ }}
\newcommand{\note}{\noindent {\bf Note:\ }}
\newcommand{\examp}{\noindent {\bf Example:\ }}
\newcommand{\Lower}[2]{\smash{\lower #1 \hbox{#2}}}
\newcommand{\ben}{\begin{enumerate}}
\newcommand{\een}{\end{enumerate}}
\newcommand{\bi}{\begin{itemize}}
\newcommand{\ei}{\end{itemize}}
\newcommand{\hp}{\hspace{.2in}}
\newtheorem{lw}{Proposition 3.1, Lo and Weng (1989)}
\newtheorem{thm}{Theorem}[section]
\newtheorem{defin}{Definition}[section]
\newtheorem{prop}{Proposition}[section]
\newtheorem{lem}{Lemma}[section]
\newtheorem{cor}{Corollary}[section]
\newcommand{\rb}[1]{\raisebox{1.5ex}[0pt]{#1}}
\newcommand{\mc}{\multicolumn}
%Mathrsfs Font
\newcommand{\Bcr}{\mathscr{B}}
\newcommand{\Ucr}{\mathscr{U}}
\newcommand{\Gcr}{\mathscr{G}}
\newcommand{\Dcr}{\mathscr{D}}
\newcommand{\CS}{\mathscr{C}}
\newcommand{\Fcr}{\mathscr{F}}
\newcommand{\Icr}{\mathscr{I}}
\newcommand{\Lcr}{\mathscr{L}}
\newcommand{\Mcr}{\mathscr{M}}
\newcommand{\Ncr}{\mathscr{N}}
\newcommand{\Pcr}{\mathscr{P}}
\newcommand{\Qcr}{\mathscr{Q}}
\newcommand{\Scr}{\mathscr{S}}
\newcommand{\Tcr}{\mathscr{T}}
\newcommand{\Xcr}{\mathscr{X}}
\newcommand{\Vcr}{\mathscr{V}}
\newcommand{\Ycr}{\mathscr{Y}}
%Mathbb Font
\newcommand{\E}{\mathbb{E}}
\newcommand{\F}{\mathbb{F}}
\newcommand{\I}{\mathbb{I}}
\newcommand{\Q}{\mathbb{Q}}
\newcommand{\X}{\mathbb{X}}
\newcommand{\Pe}{\mathbb{P}}
\newcommand{\M}{\mathbb{M}}
\newcommand{\Wbb}{\mathbb{W}}
\def\Beta{\text{Beta}}
\def\Dir{\text{Dirichlet}}
\def\DP{\text{DP}}
\def\P{{\bf p}}
\def\fhat{\widehat{f}}
\def\GA{\text{gamma}}
\def\ind{\stackrel{\mathrm{ind}}{\sim}}
\def\iid{\stackrel{\mathrm{iid}}{\sim}}
\def\J{{\bf J}}
\def\K{{\bf K}}
\def\min{\text{min}}
\def\N{\text{N}}
\def\p{{\bf p}}
\def\U{{\bf U}}
\def\W{{\bf W}}
\def\T{{\bf T}}
\def\y{{\bf y}}
\def\t{{\bf t}}
\def\m{{\bf m}}
\def\X{{\bf X}}
\def\Y{{\bf Y}}
\def\tps{\mbox{\scriptsize ${\theta H}$}}   %   smaller "\psi"-vector
\def\ups{\mbox{\scriptsize ${P_{\theta}(g)}$}}   %   smaller "\psi"-vector
\def\vps{\mbox{\scriptsize ${\theta}$}}   %   smaller "\psi"-vector
\def\vups{\mbox{\scriptsize ${\theta >0}$}}   %   smaller "\psi"-vector
\def\hps{\mbox{\scriptsize ${H}$}}   %   smaller "\psi"-vector
\def\rps{\mbox{\scriptsize ${(\theta+1/2,\theta+1/2)}$}}   %   smaller "\psi"-vector
\def\sps{\mbox{\scriptsize ${(1/2,1/2)}$}}   %   smaller "\psi"-vector
\newcommand{\reals}{{\rm I\!R}}
\newcommand{\PR}{{\rm I\!P}}
\def\Z{{\bf Z}}
\def\yy{{\mathcal Y}}
\def\rr{{\mathcal R}}
\def\BP{\text{beta}}
\def\ts{\tilde{t}}
\def\js{\tilde{J}}
\def\gs{\tilde{g}}
\def\fs{\tilde{f}}
\def\ys{\tilde{Y}}
\def\ps{\tilde{\mathcal {P}}}
\def\Report{Lancelot F. James}
\def\Author{Integrated Volatility}
\pagestyle{myheadings}
\markboth{\Author}{\Report}
\thispagestyle{empty}
\bct\Heading  A note on exact likelihoods of the Carr-Wu models\\
for leverage effects and volatility in financial
economics.\lbk\lbk\smc Lancelot F. James\footnote{ \eightit AMS
2000 subject classifications.
               \rm Primary 62G05; secondary 62F15.\\
\eightit Corresponding authors address.
                \rm The Hong Kong University of Science and Technology,
Department of Information and Systems Management, Clear Water Bay,
Kowloon, Hong Kong.
\rm lancelot\at ust.hk\\
\indent\eightit Keywords and phrases.
                \rm
          Leverage effects,
          Mixture of Normals,
          Ornstein-Uhlenbeck Process,
          Poisson Random Measure,
          Stochastic Volatility.
          }
\lbk\lbk \BigSlant The Hong Kong University of Science and
Technology\rm \lbk %(\today)%
\ect \Quote Recently Carr and Wu~(2004, 2005) and also Huang and
Wu~(2004) show that most stochastic processes used in traditional
option pricing models can be cast as special cases of time-changed
L\'evy processes. In particular these are models which can be
tailored to exhibit correlated jumps in both the log price of
assets and the instantaneous volatility. Naturally similar to a
recent work of Barndorff-Nielsen and Shephard~(2001a, b), such
models may be used in a likelihood based framework. These
likelihoods are based on the unobserved integrated volatility,
rather than the instantaneous volatility. James~(2005) establishes
general results for the likelihood and estimation of a large class
of such models which include possible leverage effects. In this
note we show that exact expressions for likelihood models based on
generalizations of Huang and Wu~(2004) and Carr and Wu~(2005),
follow essentially from the arguments in Theorem 5.1 in
James~(2005) with some slight modification. We show that that an
explicit likelihood for any of these types of models only requires
knowledge of the characteristic functional of a suitably defined
linear functional. This serves to formally verify a claim made by
James~(2005).
 \EndQuote
%\baselineskip14pt
%\begin{document}
\rm
%\newpage
\section{Introduction}
Recently Carr and Wu~(2004, 2005) and also Huang and Wu~(2004)
show that most stochastic processes used in traditional option
pricing models can be cast as special cases of time-changed L\'evy
processes. This includes for instance the models of Duffie, Pan,
and Singleton~(2000) and leverage effects model of
Barndorff-Nielsen and Shephard~(2001a, b). In particular these are
models which can be tailored to exhibit correlated jumps in both
the log price of assets and the instantaneous volatility.
Naturally similar to a recent work of Barndorff-Nielsen and
Shephard~(2001a, b), such models may be used in a likelihood based
framework. James~(2005) establishes general results for the
likelihood and estimation of a large class of such models, based
on quite general linear functionals of Poisson random measures,
which include possible leverage effects. In this note we show that
exact expressions for likelihood models based on generalizations
of Carr and Wu~(2005) and Huang and Wu~(2004), follow essentially
from the arguments in Theorem 5.1 in James~(2005) with some slight
modification. This serves to formally verify a claim made by
James~(2005). We shall be rather brief in our exposition and refer
the reader to the above mentioned works for further references and
motivation.
   Huang and Wu~(2004) and Carr and Wu~(2005) proposed a model
   for  the log price of assets which can be written as
   \Eq
   x^{*}(t)=(r-q)t+J_{1}(\tau(t))+J_{2}(\gamma(t))+\beta\tau(t)+\alpha \gamma(t)+\sigma
   W_{1}(\tau(t))+\sigma W_{2}(\gamma(t))
\label{model}
   \EndEq
   where $(J_{1}, J_{2})$ are independent pure jump L\'evy
   processes,$(W_{1},W_{2})$ are independent standard Brownian
   motion independent of $(J_{1},J_{2})$. Furthermore
   $(\tau,\gamma)$ are non-negative random time
   changes, independent of the above processes. The independence
   property can be relaxed. See Carr and Wu~(2004).
   An example of $\tau$ and $\gamma$ are the integrated Ornstein-Uhlenbeck models
   of Barndorff-Nielsen and Shephard~(2001a, b) which are used to
   model the integrated stochastic volatility. For notational
   convenience we shall hereafter set $\sigma=1$. Assuming
   conditional independence across intervals $[(i-1)\Delta,i\Delta]$
   for $i=1,\ldots, n$ and $\Delta>0$, define $\tau_{i}:=\tau((i\Delta))-\tau(((i-1)\Delta))$ and
   $\gamma_{i}:=\gamma((i\Delta))-\gamma((i-1)\Delta))$, additionally define
$$
J_{i,1}=J_{1}(\tau((i\Delta)))-J_{1}(\tau((i-1)\Delta))
$$
and
$$
J_{i,2}=J_{2}(\gamma((i\Delta)))-J_{2}(\gamma((i-1)\Delta)).
$$
It follows from~\mref{model} that one may define a likelihood
model for the aggregate returns
$X_{i}=x^{*}(i\Delta)-x^{*}((i-1)\Delta)$ which is
   conditionally Normal given $(J_{1},J_{2},\tau,\gamma)$.
   Assuming that these quantities depend on some unknown Euclidean parameter $\theta$,
   this is in essence the type of framework addressed in James~(2005)
   except for the appearance of the time-changed $(J_{1},J_{2})$.
   Before we proceed to derive an exact form of this present
   model, we will further generalize~\mref{model} by modeling the
   $\tau$ and $\gamma$ as linear functionals of a quite arbitrary
   Poisson random measure, say $N$,  on a Polish space $\Vcr$. In this way
   the $t$ in the model~\mref{model} may be replaced by a more
   abstract notion of time. Note that we can always choose $N$ on
   a big enough space so that $\tau$ and $\gamma$ are either
   independent or dependent.
   \Remark Note that we take the quite general Poisson framework
   for explicit flexible concreteness. This includes both discrete
   and continuous quantities, such as shot-noise processes on
   abstract spaces.
   As we shall see, on a more
   abstract level we simply can specify $\tau$ and $\gamma$ so
   that one knows the explicit characteristic functional of linear
   combinations of such processes. \EndRemark
\section{Exact expression for the marginal likelihood}
First as in James~(2005) let $N$ denote a Poisson random measure
on some Polish space $\Vcr$ with mean intensity,
$$
\E[N(dx)|\nu]=\nu(dx).
$$
We denote the Poisson law of $N$ with intensity $\nu$ as
$\Pe(dN|\nu)$. The Laplace functional for $N$ is defined as
$$
\E[{\mbox e}^{-N(f)}|\nu]=\int_{\M}{\mbox
e}^{-N(f)}\Pe(dN|\nu)={\mbox e}^{-\Lambda(f)}
$$
where for any positive $f$, $N(f)=\int_{\Vcr}f(x)N(dx)$ and
$\Lambda(f)=\int_\Vcr(1-{\mbox e}^{-f(x)})\nu(dx).$ $\M$ denotes
the space of boundedly finite measures on $\Vcr.$
 We suppose that
$\tau=N(h_{i})$, for $i=1,\ldots, n$ where $h_{1},\ldots,h_{n}$
are positive measureable functions on $\Vcr$. Similarly we suppose
that for each $i$ $\gamma_{i}=N(g_{i})$ where $g_{i}$ are positive
measureable functions on $\Vcr$. With this specification it
follows from~\mref{model} that the generalized notion of aggregate
returns
$X_{i}|\tau_{i},\gamma_{i},J_{i,1},J_{i,2},\theta,\alpha,\beta,r,q$
are conditionally independent Normal random variables expressible
as \Eq X_{i}=(r-q)\Delta
+J_{i,1}+J_{i,2}+\beta(\tau_{i}+\gamma_{i})+(\alpha-\beta)\gamma_{i}+\sqrt{\tau_{i}+\gamma_{i}}\epsilon_{i}
\label{cond}\EndEq where $\epsilon_{i}$ are independent standard
Normal random variables. Now for $f_{i}(x):=h_{i}(x)+g_{i}(x)$ on
$\Vcr$, set
$$\tau^{*}_{i}:=\tau_{i}+\gamma_{i}=N(f_{i})$$
and set $\mu=(r-q)$ and $\alpha-\beta=\rho$. Then the conditional
Normal density of each $X_{i}$, say  $\phi(X_{i}|\mu\Delta
+J_{i,1}+J_{i,2}+\rho\gamma_{i}+\beta\tau^{*}_{i},\tau^{*}_{i})$,
can be written as \Eq {\mbox
e}^{(A_{i}-J_{i,1}-J_{i,2}-\rho\gamma_{i})\beta
}\frac{1}{\sqrt{2\pi}}\frac{1}{\sqrt{\tau^{*}_{i}}}[{\mbox
e}^{-{(A_{i}-J_{i,1}-J_{i,2}-\rho
\gamma_{i})}^{2}/(2\tau^{*}_{i})}]{\mbox
e}^{-\tau^{*}_{i}\beta^{2}/2} \label{den}\EndEq where
$A_{i}=X_{i}-\mu\Delta$. It is not difficult to see that if
$J_{i,1}$ and $J_{i,2}$ were removed, the subsequent likelihood
model is a special case of the models handled by James~(2005,
Theorem 5.1). James~(2005) remarks that the models such
as~\mref{den} pose no additional difficulties as long as one knows
the characteristic functional of say $J_{i,1}$ and $J_{i,2}$ and
their sums. However $J_{1}$ and $J_{2}$ are pure jump L\'evy
processes and hence by way of the known L\'evy -Khinchine formula
and independence relative to $\tau$ and $\gamma$, this point is a
trivial matter. In particular for any real or complex number,
$\omega$ write for j=$1,2$,
$$
\E[{\mbox e}^{-\omega [J_{j}(t)-J_{j}(s)]}]={\mbox
e}^{-(t-s)\psi_{j}(\omega)}
$$
where the explicit form of $\psi_{j}(\omega)$ is given by the
L\'evy-Khinchine formula which can be found in Carr and Wu~(2005),
but is otherwise similar to the function $\Lambda$. It then
follows that for each $i$
$$
\E[{\mbox e}^{-\omega J_{i,1}}]=\E[{\mbox
e}^{-\tau_{i}\psi_{1}(\omega)}]={\mbox
e}^{-\Lambda(h_{i}\psi_{1}(\omega))} {\mbox { and }}\E[{\mbox
e}^{-\omega J_{i,2}}]={\mbox
e}^{-\Lambda(g_{i}\psi_{2}(\omega))}$$ Note also conditional on
$N$, for possibly complex valued numbers
$(\omega_{1},\ldots,\omega_{n})$, \Eq \E[\prod_{i=1}^{n}{\mbox
e}^{-\omega_{i} J_{i,1}}|N]=\prod_{i=1}^{n}{\mbox
e}^{-\tau_{i}\psi_{1}(\omega_{i})}{\mbox { and
}}\E[\prod_{i=1}^{n}{\mbox e}^{-\omega_{i}
J_{i,2}}|N]=\prod_{i=1}^{n}{\mbox
e}^{-\gamma_{i}\psi_{2}(\omega_{i})}. \label{key}\EndEq How we
shall proceed is to first evaluate everything conditional on $N$.
Our results will then boil down to expectation of exponential sum
of terms of the form, \Eq \tau_{i}[\psi_{1}(\beta+\xi
y_{i})+(\beta^{2}+y^{2}_{i})/2]+ \gamma_{i}[\psi_{2}(\beta+\xi
y_{i})+\rho(\beta+\xi y_{i})+(\beta^{2}+y^{2}_{i})/2]
\label{done}\EndEq where $\xi$ is the imaginary number. This is
similar to the case of James~(2005). Now for real valued numbers
$(y_{1},\ldots,y_{n})$, set
$$
\Omega_{n}(x)=\sum_{i=1}^{n}[\psi_{1}(\beta+\xi
y_{i})+(\beta^{2}+y^{2}_{i})/2]h_{i}(x)
$$
and
$$
\Upsilon_{n}(x)=\sum_{i=1}^{n}[\psi_{2}(\beta+\xi
y_{i})+\rho(\beta+\xi y_{i})+(\beta^{2}+y^{2}_{i})/2]g_{i}(x)
$$
Note that the sum over the terms in ~\mref{done} is equivalent in
distribution to $N(\Omega_{n}+\Upsilon_{n}).$ We now state the
form of the likelihood.
\begin{thm} Suppose that $N$ is a Poisson random measure with intensity $\nu$ on
$\Vcr$. Furthermore suppose that $\tau_{i}$ and $\gamma_{i}$,
defined above, are chosen such that
$\Lambda(\Omega_{n}+\Upsilon_{n})<\infty.$ Then the joint marginal
density or likelihood of
$X_{1},\ldots,X_{n}|\mu,\beta,\theta,\rho$, determined
by~\mref{cond} and~\mref{den} is given by,
$$\Lcr(\X|\mu,\beta,\theta,\rho)= \frac{{\mbox e}^{n{\bar A}\beta
}}{{(2\pi)}^{n}}\int_{{\mathbb R}^{n}} {\mbox
e}^{-\Lambda(\Omega_{n}+\Upsilon_{n})} \prod_{i=1}^{n}{\mbox
e}^{\xi A_{i}y_{i}}dy_{i}.$$ Where ${\bar
A}=\sum_{i=1}^{n}A_{i}/n$. The result applies for the case where
$N$ is not necessarily Poisson, by replacing ${\mbox
e}^{-\Lambda(\Omega_{n}+\Upsilon_{n})}$ with $\E[{\mbox
e}^{-N(\Omega_{n}+\Upsilon_{n})}].$ \qed\end{thm} \Proof The proof
of this result is simply a slight variation of Theorem 5.1 in
James~(2005). For completeness we give many of the same details.
Here we use the fact that for each $i$ one has the identity
deduced from the characteristic function of a Normal distribution,
with mean 0 and variance $1/\tau^{*}_{i}$, evaluated at
$\varpi_{i}=A_{i}-J_{i,1}-J_{i,2}-\rho \gamma_{i}$. That is,
$$
\frac{1}{\sqrt{2\pi}}\int_{-\infty}^{\infty}{\mbox
e}^{\xi\varpi_{i}y_{i}-\tau^{*}_{i}y^{2}_{i}/2}dy_{i}=\frac{1}{\sqrt{\tau^{*}_{i}}}{\mbox
e}^{-{(\varpi_{i})}^{2}/2\tau^{*}_{i}}
$$
Now the result proceeds by substituting this expression
in~\mref{den} and applying Fubini's theorem . One then integrates
out the expressions involving $(J_{i,1}, J_{i,2})$ conditionally
on $N$, which results in using~\mref{key}. The rest now is
precisely as the proof of  Theorem 5.1 in James~(2005).  That is
after rearranging terms it remains to calculate the expectation of
${\mbox e}^{-N(\Omega_{n}+\Upsilon_{n})}$ \EndProof \Remark Note
that these arguments may be used to directly evaluate the density
of $x^{*}(t)$ given observations $X_{1},\ldots,X_{n}$, which might
be interesting in an option pricing context. Statistical
estimation follows along the lines of James~(2005). \EndRemark
 \vskip0.2in \centerline{\Heading References}
\vskip0.2in \tenrm
\def\smc{\tensmc}
\def\sl{\tensl}
\def\bf{\tenbold}
\baselineskip0.15in
\Ref \by Barndorff-Nielsen, O.E. and Shephard, N. \yr 2001a \paper
Ornstein-Uhlenbeck-based models and some of their uses in
financial economics \jour \JRSSB \vol 63 \pages 167-241 \EndRef
\Ref \by Barndorff-Nielsen, O.E. and Shephard, N. \yr 2001b \paper
Modelling by L\'evy processes for financial econometrics. In
L\'evy processes. Theory and applications. Edited by Ole E.
Barndorff-Nielsen, Thomas Mikosch and Sidney I. Resnick. p.
283-318. Birkh\"auser Boston, Inc., Boston, MA \EndRef \Ref \by
Carr, P. and Wu, L. \yr 2004 \paper Time-changed L\'evy processes
and option pricing \jour Journal of Financial Economics \vol 71
\pages 113-141 \EndRef
\Ref \by Carr, P. and Wu, L. \yr 2005 \paper Stochastic skew in
currency options. Manuscript available at
http://faculty.baruch.cuny.edu/lwu/ \EndRef
 \Ref \by Duffie, D., Pan, J. and Singleton, K.\yr 2000 \paper Transform analysis and
asset pricing for affine jump diffusions \jour Econometrica \vol
68 \pages 1343-1376 \EndRef
\Ref \by Huang, J.Z. and Wu, L. \yr 2004 \paper Specification
analysis of option pricing models based on time-changed L\'evy
processes \jour Journal of Finance \vol 59 \pages 1405-1439
\EndRef \Ref\by James, L.F. \yr 2005 \paper Analysis of a class of
likelihood based continuous time stochastic volatility models
including Ornstein-Uhlenbeck models in financial economics.
arXiv:math.ST/0503055 \EndRef
\medskip
\smc
\Tabular{ll}
Lancelot F. James\\
The Hong Kong University of Science and Technology\\
Department of Information and Systems Management\\
Clear Water Bay, Kowloon\\
Hong Kong\\
\rm lancelot\at ust.hk\\
\EndTabular
\end{document}